\documentclass[leqno,11pt]{article}  
\usepackage{amsmath}
\usepackage{amssymb}

\usepackage[german,english]{babel}
\usepackage{amsthm}

\title{A note on multivariable $(\varphi,\Gamma)$-modules}

\author{\textsc{Elmar Grosse-Kl\"onne}}
\date{}

\theoremstyle{plain} 
\newtheorem{satz}{Theorem}
\newtheorem{kor}[satz]{Corollary}  
\newtheorem{lem}[satz]{Lemma}  
\newtheorem{pro}[satz]{Proposition}  
 
\theoremstyle{remark}

\theoremstyle{definition}

\newcommand{\0}{\ensuremath{\overrightarrow{0}}}

\begin{document}

\maketitle

%


\begin{abstract} Let $F/{\mathbb Q}_p$ be a finite field extension, let $k$ be a field of characteristic $p$. Fix a Lubin Tate group $\Phi$ for $F$ and let $\Gamma_{\bullet}=\Gamma\times\cdots\times\Gamma$ with $\Gamma={\mathcal O}_F^{\times}$ act on $k[[t_1,\ldots,t_n]][\prod_it_i^{-1}]$ by letting $\gamma_i$ (in the $i$-th factor $\Gamma$) act on $t_i$ by insertion of $t_i$ into the power series attached to $\gamma_i$ by $\Phi$. We show that $k[[t_1,\ldots,t_n]][\prod_it_i^{-1}]$ admits no non-trivial ideal stable under $\Gamma_{\bullet}$, thereby generalizing a result of Z\'{a}br\'{a}di (who had treated the case where $\Phi$ is the multiplicative group). We then discuss applications to $(\varphi,\Gamma)$-modules over $k[[t_1,\ldots,t_n]][\prod_it_i^{-1}]$.\\ 
\end{abstract}

\begin{center}{\Large{\bf Introduction}}\end{center} 

In recent years, Fontaine's by now classical theory of
$(\varphi,\Gamma)$-modules classifying $p$-adic and $p$-torsion
representations of the absolute Galois group of a local field with residue
characteristic $p$ has been enriched by several new ideas and
variations. Among them are various concepts involving
$(\varphi,\Gamma)$-modules in several free variables (as opposed to the
classical use of
a single variable). One such concept is due to Z\'{a}br\'{a}di
\cite{zabr}. Motivated by constructions related to the $p$-adic local
Langlands program he introduced $(\varphi,\Gamma)$-modules over rings
$A[[t_1,\ldots,t_n]][\prod_it_i^{-1}]$ where $A$ is a finite quotient of the
ring of integers in a finite extension of ${\mathbb Q}_p$ and where for each
$1\le i\le n$ we have an action of $\Gamma\cong{\mathbb Z}_p^{\times}$
and $\varphi$ given by inserting $t_i$ into the formal group law ${\mathbb
  G}_m$. Later he showed \cite{zabgal} that these $n$-variable
$(\varphi,\Gamma)$-modules classify $p$-adic and $p$-torsion
representations of the $n$-fold power of ${\rm Gal}(\overline{\mathbb
  Q}_p/{\mathbb Q}_p)$. Let
$k$ be the residue field of $A$. The fact that
$k[[t_1,\ldots,t_n]][\prod_it_i^{-1}]$ is no longer a field if $n>1$ a priori
causes numerous technical problems when trying to generalize classical methods
to this context. However, Z\'{a}br\'{a}di overcomes them all at once with a
clever argument: He shows that $k[[t_1,\ldots,t_n]][\prod_it_i^{-1}]$ contains
no non-trivial ideal simultaneously stable under all the $\Gamma$-actions.

In another direction, Fontaine's classical theory of
$(\varphi,\Gamma)$-modules has been generalized by Kisin, Ren and Schneider \cite{kisren}
\cite{peterlec}, see also \cite{berfou}, \cite{fouxie}, \cite{schven}, by
substituting the formal group law ${\mathbb
  G}_m$ with more general Lubin Tate group laws. 

One
purpose of this note is to explain how the aforementioned strategy of
Z\'{a}br\'{a}di can be implemented in the context which commonly
generalizes both his one as well as that of Kisin, Ren and
Schneider. This will allow, for example, generalizations of results from
\cite{zabr} and notably \cite{zabgal} from ${\mathbb Q}_p$ to finite
extensions of ${\mathbb Q}_p$.

The second purpose is to introduce resp. clarify a construction of $(\varphi,\Gamma)$-modules over
$A[[t_1,\ldots,t_n]][\prod_it_i^{-1}]$ (with respect to a general a Lubin Tate
group law) out of certain torsion
$A[[t_1,\ldots,t_n]]$-modules endowed with actions by
$\Gamma$ and $\varphi$. The latter are often naturally assigned to smooth
representations of $p$-adic reductive groups, see \cite{breuil},
\cite{zabr}, \cite{fungl} and further forthcoming work of the author. The
construction rests on employing $\psi$-operators on
$A[[t_1,\ldots,t_n]][\prod_it_i^{-1}]$ which so far in the literature has been done
only in the classical ${\mathbb G}_m$-context.

{\it Acknowledgements:} I thank Gergely Z\'{a}br\'{a}di for his careful
reading of an earlier draft of the proof of Theorem \ref{gergelytrick} and for
further discussions on the topic. I thank the referees for their valuable comments.\\ 

{\bf Notation:} Let $K/F/{\mathbb Q}_p$ be finite field extensions. Let $q$ be the number of
elements of the
residue field of $F$, let $\pi$ be a uniformizer in the
ring of integers ${\mathcal O}_F$ of $F$. Let $A$ be a finite quotient of the ring ${\mathcal O}_K$ of integers of $K$
and let $k$ be its residue field.

Put $\Gamma={\mathcal O}_F^{\times}$. There is a unique Lubin-Tate group for $F$ with respect to $\pi$; fixing a coordinate $t$ we write $\Phi(t)$ for the corresponding Lubin-Tate formal power series. For $\gamma\in\Gamma$ let
$[\gamma]_{\Phi}(t)\in{\mathcal O}_F[[t]]$ denote the power series describing
the action of $\gamma$ in the Lubin-Tate group. Let $D$ be a finite set. Put $${A}[[t_\bullet]]={A}[[t_d]]_{d\in
  D},\quad\quad {A}((t_{\bullet}))={A}[[t_{\bullet}]][\prod_{d\in
    D}t_d^{-1}].$$For each $d\in D$ let $\Gamma_d$ be a copy of $\Gamma$. For
$\gamma\in\Gamma$ let $\gamma_d$ denote the element in $\Gamma_{\bullet}=\prod_{d\in
  D}\Gamma_d$ whose $d$-component is $\gamma$ and whose other components are trivial. The formulae $$\gamma_{d}*
t_{d}=[\gamma]_{\Phi}(t_{d})\quad\quad\quad \gamma_{d_1}*
t_{d_2}=t_{d_2}$$with $\gamma\in\Gamma$ and $d, d_1, d_2\in D$ such that
$d_1\ne d_2$ define an action of $\Gamma_{\bullet}$ on ${A}((t_{\bullet}))$. It induces an action of $\Gamma_{\bullet}$ on $k((t_{\bullet}))$.\\

The following theorem was proven in \cite{zabr} in the case $\Phi(t)=(1+t)^p-1$ (and hence $F={\mathbb Q}_p$). That proof made critical use of the fact
that for $\Phi(t)=(1+t)^p-1$ the coefficients of the power series
$[\gamma]_{\Phi}(t)$ are explicitly known (by a closed formula), and more
importantly, that $[\gamma]_{\Phi}(t)$ is in fact a polynomial for
sufficiently many $\gamma\in\Gamma$ (namely the $\gamma=1+p^n$ with
$n>0$). That proof breaks down for other $\Phi$'s.

\begin{satz}\label{gergelytrick} There exists no
  nontrivial ideal in $k((t_{\bullet}))$ stable under the action of $\Gamma_{\bullet}$.
\end{satz}

{\sc Proof:} Throughout we fix an element $\gamma\in 1+\pi{\mathcal O}_F$ with $\gamma\ne 1$. We then have
$$[\gamma]_{\Phi}(t)-t+ut^{n}\in t^{n+1} k[[t]]$$ for some
$n\ge 2$, some $u\in k^{\times}$. For $d\in D$ and $h\in k[[t_{\bullet}]]-\{0\}$ let ${\rm ord}_{d}(h)$ be the
maximal integer $o$ with $$h\in t_{d}^ok[[t_{\bullet}]].$$For $m\ge0$ and $f\in
k[[t_{d}]]$ with ${\rm
  ord}_{d}(f)>0$ and $h\in k[[t_{d'}]]_{d'\ne
    d}-\{0\}$ we have $$\gamma_{d}*(t_{d}^mhf)-t_{d}^mhf+u(t_{d}^mhf)^{n}\in
(t_{d}^mhf)^{n+1}k[[t_{\bullet}]],$$i.e. $${\rm
  ord}_{d}(\gamma_{d}*(t_{d}^mhf)-t_{d}^mhf)=n({\rm ord}_{d}(f)+m),$$i.e.\begin{gather}{\rm
  ord}_{d}(t_{d}^{-{n}{\rm ord}_{d}(f)}(\gamma_{d}*(t_{d}^mhf)-
t_{d}^mhf))=n m.\label{dochkri}\end{gather}On the other hand, from $$t_{d}^{-{n}{\rm ord}_{d}(f)}(\gamma_{d}*(t_{d}^mhf)-
t_{d}^mhf)=h t_{d}^{-{n}{\rm ord}_{d}(f)}(\gamma_{d}*(t_{d}^mf)-
t_{d}^mf)$$we deduce\begin{gather}{\rm
  ord}_{d'}(t_{d}^{-{n}{\rm ord}_{d}(f)}(\gamma_{d}*(t_{d}^mhf)-
t_{d}^mhf))={\rm
  ord}_{d'}(h)\quad\quad\mbox{ for }d'\ne d.\label{dochkri1}\end{gather}Let $g\in
k[[t_{\bullet}]]^{\times}$. We ask for the values$${\rm
  ord}_{d'}(t_{d}^{-{n}{\rm ord}_{d}(f)}(\gamma_{d}*(t_{d}^mhf)-
t_{d}^mhf))$$where $t_{d}^mh$ runs through all the monomials
appearing in $g$. It follows from formulae (\ref{dochkri}) and
(\ref{dochkri1}) that for the constant monomial this value is $=0$ for all
$d'\in D$, but that
for all the other monomials there is some $d'\in D$ for which this value is $>0$. We deduce\begin{gather}t_{d}^{-{n}{\rm ord}_{d}(f)}(\gamma_{d}*(
gf)-gf)\in
k[[t_{\bullet}]]^{\times}\label{karsa}\\\mbox{ for }g\in
k[[t_{\bullet}]]^{\times}\mbox{ and } f\in
k[[t_{d}]]\mbox{ with }{\rm
  ord}_{d}(f)>0.\notag\end{gather}We may write elements in $k((t_{\bullet}))$
as \begin{gather}h=g+\sum_{j\in J}f_j'\prod_{d\in D_j}f_{d,j}\label{tensorsum}\\\mbox{
    with } g\in
  k[[t_{\bullet}]], D_j\subset D, f_j'\in
  k[[t_{d'}]]_{d'\in D-D_j}, f_{d,j}\in
  k((t_{d}))^{\times}\notag\end{gather}for finite index sets $J$. Let
$I\subset k((t_{\bullet}))$ be a non zero ideal, stable under $\Gamma_{\bullet}$. Multiplying an arbitrary
non zero element in $I$ with a suitable
monomial we see that $I$ contains an element $h$ written in the form
(\ref{tensorsum}) such
that $g$ is a unit in $k[[t_{\bullet}]]$. Among all these elements $h$ and
all expressions (\ref{tensorsum}) for such $h$ we choose one for which $J$
has minimal cardinality. We claim that $J$ is in fact empty. Assuming the
contrary, we find some $j_0\in J$ and some $d_0\in D_{j_0}$ such that $f_{d_0, j_0}\notin
k[[t_{d_0}]]$. (If there was no such $j_0$ then we would have $h\in
k[[t_{\bullet}]]$, hence even $h\in
k[[t_{\bullet}]]^{\times}$ after possibly multiplying with a monomial in
$k((t_{\bullet}))$. But this would mean that $J$ is empty.) Notice that this implies $f_{d_0,j_0}^{-1}\in
k[[t_{d_0}]]$ and ${\rm
  ord}_{d_0}(f_{d_0, j_0}^{-1})>0$. We then consider the
element \begin{gather}t_{d_0}^{-{n}{\rm
      ord}_{d_0}(f_{d_0,j_0}^{-1})}(\gamma_{d_0}*(f_{d_0,j_0}^{-1}h)-f_{d_0,j_0}^{-1}h)\label{newexempl}\end{gather}
of $I$. It is the sum of \begin{gather}t_{d_0}^{-{n}{\rm
      ord}_{d_0}(f_{d_0,j_0}^{-1})}(\gamma_{d_0}*(f_{d_0,j_0}^{-1}g)-f_{d_0,j_0}^{-1}g)\label{exempl}\end{gather}and
of all\begin{gather}t_{d_0}^{-{n}{\rm ord}_{d_0}(f_{d_0,j_0}^{-1})}(\gamma_{d_0}*(f_{d_0, j_0}^{-1}f_j'\prod_{d\in D_j}f_{d,j})-f_{d_0,j_0}^{-1}f_j'\prod_{d\in D_j}f_{d,j})\label{nonexempl}\end{gather}for ${j\in J}$. By (\ref{karsa}), the
element (\ref{exempl}) is a unit in $k[[t_{\bullet}]]$. Each
  element (\ref{nonexempl}) may be written
  as$$f_j' t_{d_0}^{-{n}{\rm
      ord}_{d_0}(f_{d_0, j_0}^{-1})}(\gamma_{d_0}*(f_{d_0,j_0}^{-1}f_{d_0,j})-f_{d_0, j_0}^{-1}f_{d_0,j})\prod_{d\in
    D_j-\{d_0\}}f_{d,j}.$$The factor
  $$t_{d_0}^{-{n}{\rm
      ord}_{d_0}(f_{d_0,j_0}^{-1})}(\gamma_{d_0}*(f_{d_0,j_0}^{-1}f_{d_0,j})-f_{d_0,j_0}^{-1}f_{d_0,j})$$
  belongs to $k((t_{d_0}))$, and for $j=j_0$ it even vanishes. This shows that
  $J$ was not chosen minimally, refuting our assumption. Thus, $J$ is
  empty. But then $h=g$ is a unit in $k[[t_{\bullet}]]$, in particular in
  $k((t_{\bullet}))$, hence $I=k((t_{\bullet}))$.\hfill$\Box$\\

 {\bf Remarks:} It is straightforward to streamline the above proof to show more: Fix $d'\in D$ and for $d\in D-\{d'\}$ let $\widetilde{\Gamma}_{d}\subset \Gamma_{d}$ be a subgroup such that $[\gamma]_{\Phi}(t)$ is non-linear for at least one $\gamma_{d}\in \widetilde{\Gamma}_{d}$; then there exists no
  nontrivial ideal in $k((t_{\bullet}))$ stable under
  $\prod_{d\in D-\{d'\}}\widetilde{\Gamma}_{d}$.\\

Consider the $A$-algebra$${A}[[t_{\bullet}]][\varphi_{\bullet},\Gamma_{\bullet}]={A}[[t_{d}]]_{d\in D}[\varphi_{d},\Gamma_{d}]_{d\in D}$$ with
commutation rules given by$$x_{d_1}\cdot y_{d_2}=y_{d_2}\cdot
x_{d_1},$$$$\gamma_{d}\cdot
\varphi_{d}=\varphi_{d}\cdot\gamma_{d},\quad\quad \gamma_{d}\cdot
t_{d}=(\gamma_{d}*t_{d})\cdot\gamma_{d},\quad\quad\varphi_{d}\cdot
t_{d}=\Phi(t_{d})\cdot \varphi_{d}$$for $\gamma\in\Gamma$ and $x,y\in\Gamma\cup\{\varphi,t\}$
and $d, d_1, d_2\in D$ with $d_1\ne d_2$. Similarly we define the $A$-algebra
${A}((t_{\bullet}))[\varphi_{\bullet},\Gamma_{\bullet}]$ and its subalgebra ${A}((t_{\bullet}))[\Gamma_{\bullet}]$.

\begin{kor}\label{herzmariens} Let ${\bf
  D}$ be an ${A}((t_{\bullet}))[\Gamma_{\bullet}]$-module, let ${\bf
  D}'\subset{\bf
  D}$ be an ${A}((t_{\bullet}))[\Gamma_{\bullet}]$-sub module, let $\alpha\in
  {A}((t_{\bullet}))$ with $\alpha{\bf
  D}'=0$. If $\alpha$ maps to a non-zero element in
  $k((t_{\bullet}))$, then ${\bf
  D}'=0$.
\end{kor}

{\sc Proof:} Let ${\mathfrak m}_A$ be the maximal ideal in $A$ and let $m\ge
0$ be maximal with the property ${\bf
  D}'\subset {\mathfrak m}_A^m{\bf
  D}$. Consider the annihilator in
${k}((t_{\bullet}))$ of the image of ${\bf
  D}'$ in ${\mathfrak m}_A^m{\bf
  D}/{\mathfrak m}_A^{m+1}{\bf
  D}$. As it contains the image of $\alpha$, it must be all of
${k}((t_{\bullet}))$ by Theorem
\ref{gergelytrick}, hence ${\bf
  D}'=0$.\hfill$\Box$\\

{\bf Definition:} A $\psi$-operator on ${A}[[t_{\bullet}]]$ is a system
$\psi_{\bullet}$ of additive maps $$\psi_{d}:{A}[[t_{\bullet}]]\longrightarrow
{A}[[t_{\bullet}]]$$ for $d\in D$ such that
$\psi_{d}(\gamma_{d'}*t_{d}))=\gamma_{d'}*(\psi_{d}(t_{d}))$
for all $\gamma\in\Gamma$ and $d, d'\in D$, such that $\psi_{d}(t_{d'})=t_{d'}$ for $d\ne
d'$ and such
that the following holds true: If we view the $\varphi_{d}$ as acting on
${A}[[t_{\bullet}]]$, then $\psi_{d}\circ\varphi_{d'}=\varphi_{d'}\circ\psi_{d}$ for $d\ne
d'$, but $$\psi_{d}(\varphi_{d}(a)x)=a\psi_{d}(x)$$ for $a,x\in {A}[[t_{\bullet}]]$.\footnote{Notice that we do not require $\psi_{d}(1)=1$.} 

\begin{lem}\label{erleicht} There is a $\psi$-operator on ${A}[[t_{\bullet}]]$ such that each $\psi_{d}$ is surjective.  
\end{lem}

{\sc Proof:} We follow a construction explained in \cite{schven} section
3, see also \cite{fungl}. We assume $K=F$ since in the general case the construction given below carries over via base extension ${\mathcal O}_F\to{\mathcal O}_K$. Let us for the moment concentrate on the case $|D|=1$ (and omit indices $(.)_{d}$). The formula $\varphi\cdot
t=\Phi(t)$ defines $\varphi$ as an injective endomorphism of ${\mathcal O}_F[[t]]$. The map$$\varphi({\mathcal O}_F[[t]])^q\longrightarrow {\mathcal O}_F[[t]],\quad (a_0,\ldots,a_{q-1})\mapsto\sum_{i=0}^{q-1}a_it^{i}$$is surjective: This can be checked modulo $\pi$, hence follows from $\Phi(t)\equiv t^q$ modulo $\pi$. The map is also injective, as follows from Proposition 1.7.3 in \cite{peterlec}. It follows that
$\psi=\frac{1}{\pi}\varphi^{-1}{\rm tr}_{{\mathcal O}_F[[t]]/\varphi({\mathcal
    O}_F[[t]])}$ defines an operator on ${\mathcal O}_F[[t]]$ satisfying
$\psi(\varphi(a)x)=a\psi(x)$. To see the commutation with the
$\Gamma$-action we proceed similarly as in \cite{schven} Remark 3.2 iv. Let ${\mathcal Z}$ denote the set of $\pi$-torsion points (in the maximal ideal of ${\mathcal O}_F$) for $\Phi$. Let $F_1$ denote the extension of $F$ generated by the elements of ${\mathcal Z}$. For $z\in{\mathcal Z}$ we have the ${\mathcal O}_F$-algebra morphism $\sigma_z:{\mathcal O}_F[[t]]\to {\mathcal O}_{F_1}[[t]]$ with $t\mapsto z+_{\Phi}t$ (where $z+_{\Phi}t$ indicates addition with respect to the formal group law $\Phi$). It follows from \cite{schven} formula (10) that ${\rm tr}_{{\mathcal O}_F[[t]]/\varphi({\mathcal
    O}_F[[t]])}=\sum_{z\in{\mathcal Z}}\sigma_z$. For $\gamma\in\Gamma$ and
$a=a(t)\in {\mathcal O}_F[[t]]$ we thus compute$$\varphi(\gamma\cdot
(\psi(a(t))))=\varphi(\psi(a([\gamma]_{\Phi}(t))))=\psi(a([\gamma]_{\Phi}([\pi]_{\Phi}(t))))$$$$=\psi(a([\pi]_{\Phi}([\gamma]_{\Phi}(t))))=\varphi(\psi(a))([\gamma]_{\Phi}(t))=\frac{1}{\pi}\sum_{z\in{\mathcal
    Z}}(\sigma_z(a))([\gamma]_{\Phi}(t))$$$$=\frac{1}{\pi}\sum_{z\in{\mathcal
    Z}}\sigma_{[\gamma^{-1}]_{\Phi}(z)}(a([\gamma]_{\Phi}(t)))=\frac{1}{\pi}\sum_{z\in{\mathcal
    Z}}\sigma_{z}(a([\gamma]_{\Phi}(t)))=\varphi(\psi(\gamma\cdot
a(t))),$$hence $\gamma\cdot (\psi(a))=\psi(\gamma\cdot a)$ as $\varphi$ is
injective. To see surjectivity of $\psi$ we may assume $\Phi(t)=\pi t+t^q$ as
well as $A=k$; it is then enough to prove the following formulae (\ref{psineqp}) and (\ref{psieqqp}). If $F\ne{\mathbb Q}_p$ then for $m\in{\mathbb Z}_{\ge0}$ and
$0\le i\le q-1$ we have\begin{gather}\psi(t^{mq+i})=\left\{\begin{array}{l@{\quad:\quad}l}0
&  0\le i\le q-2\\t^m& i=q-1\end{array}\right.\label{psineqp}.\end{gather}If $F={\mathbb Q}_p$ then for $m\in{\mathbb Z}_{\ge0}$ and
$0\le i\le q-1$ we have\begin{gather}\psi(t^{mq+i})=\left\{\begin{array}{l@{\quad:\quad}l}\frac{q}{\pi}t^m& i=0\\0
&  1\le i\le q-2\\t^m& i=q-1\end{array}\right.\label{psieqqp}.\end{gather}It
is enough to
prove these formulae for $m=0$ since
$\psi(t^qa)=t\psi(a)$ for all $a\in k[[t]]$. We compute ${\rm tr}_{{\mathcal O}_F[[t]]/\varphi({\mathcal
    O}_F[[t]])}(t^{i})$ by looking at the matrix of $t^{i}$ with respect to the basis $1,t,\ldots,t^{q-1}$. Namely, for $0\le i, j\le q-1$
let us write\begin{gather}t^{i}t^j=\left\{\begin{array}{l@{\quad:\quad}l}t^{i+j}& 0\le i+j\le q-1\\(t^q+\pi
t)t^{i+j-q}- \pi t^{i+j-q+1}& q\le i+j\le 2q-2\end{array}.\right.\end{gather} For $1\le i\le q-2$ none of the $t^{i}t^j$ contributes to $\frac{1}{\pi}\varphi^{-1}{\rm tr}_{{\mathcal O}_F[[t]]/\varphi({\mathcal
    O}_F[[t]])}(t^i)$ modulo $(\pi)$. For $i=0$ we have $q$ many
summands $\frac{1}{\pi}$; their sum disappears modulo $\pi$
if and only if $F\ne{\mathbb Q}_p$. For $i=q-1$ the first line (in the above case distinction) contributes once,
and the second summand in the second line contributes $(q-1)$ times, and the outcome is as stated. Formulae (\ref{psineqp}) and (\ref{psieqqp}) are proven and the case $|D|=1$ has been settled. 

Now consider the case of a general $D$. Fix $d\in D$. We observed that ${\mathcal O}_F[[t_{d}]]$ is free of
rank $q$ over $\varphi_{d}({\mathcal O}_F[[t_{d}]])$. As ${\mathcal
  O}_F[[t_{\bullet}]]$ is the push out of the inclusions
$$\varphi_{d}({\mathcal O}_F[[t_{d}]])\longrightarrow {\mathcal
  O}_F[[t_{d}]]\quad\mbox{ and }\quad\varphi_{d}({\mathcal
  O}_F[[t_{d}]])\longrightarrow\varphi_{d}({\mathcal
  O}_F[[t_{\bullet}]])$$we may consider the operator $\frac{1}{\pi}(\varphi_{d})^{-1}{\rm tr}_{{\mathcal
    O}_F[[t_{\bullet}]]/\varphi_{d}({\mathcal
    O}_F[[t_{\bullet}]])}$ on
${\mathcal O}_F[[t_{\bullet}]]$ and proceed in the same way as we did before.\hfill$\Box$\\

An \'{e}tale $\varphi_{\bullet}$-module over
${A}((t_{\bullet}))$ is an ${A}((t_{\bullet}))[\varphi_{\bullet}]$-module ${\bf D}$
which is finitely generated over ${A}((t_{\bullet}))$ such that for each $d\in D$ the linearized structure map$${\rm id}\otimes\varphi_{d}:{A}((t_{\bullet}))\otimes_{\varphi_{d},{A}((t_{\bullet}))} {\bf
  D}\longrightarrow{\bf
  D}$$is bijective. An \'{e}tale $(\varphi_{\bullet},\Gamma_{\bullet})$-module over
${A}((t_{\bullet}))$ is an
${A}((t_{\bullet}))[\varphi_{\bullet},\Gamma_{\bullet}]$-module whose
underlying $\varphi_{\bullet}$-module is \'{e}tale. In the case $|D|=1$ we
drop the indices $(.)_{d}$ resp. $(.)_{\bullet}$ and simply talk about \'{e}tale $(\varphi,\Gamma)$-modules over
${A}((t))$.\\

{\bf Remark:} The action of $\Gamma$ on an \'{e}tale $(\varphi_{\bullet},\Gamma_{\bullet})$-module is automatically continuous for the weak topology, cf. Theorem 2.2.8 of \cite{peterlec} (for $|D|=1$).

\begin{lem} The category of \'{e}tale
$(\varphi_{\bullet},\Gamma_{\bullet})$-modules over ${A}((t_{\bullet}))$ is abelian.
\end{lem}

{\sc Proof:} With Theorem \ref{gergelytrick} available, the proof given in \cite{zabr} Proposition 2.5 (which deals with the case $\Phi(t)=(1+t)^p-1$) can be literally copied over.\hfill$\Box$\\

For an ${A}$-module $\Delta$ we write $\Delta^*={\rm
  Hom}_{A}(\Delta,{A})$. An ${A}[[t_{\bullet}]]$-module $\Delta$ is called admissible if
$$\Delta[t_{\bullet}]=\{x\in\Delta\,;\,t_{d}x=0\mbox{ for each }d\in D\}$$
is a finitely generated ${A}$-module.

We may and do assume $\Phi(t)=\pi t+t^q$. We fix the $\psi$-operator on ${A}[[t_{\bullet}]]$ constructed in the proof of
Lemma \ref{erleicht}; thus each $\psi_{d}$ satisfies formula (\ref{psineqp}) resp. (\ref{psieqqp}).

\begin{pro}\label{nopsi} Let $\Delta$ be a finitely
  generated ${A}[[t_{\bullet}]][\varphi_{\bullet},\Gamma_{\bullet}]$-module which is admissible over
  ${A}[[t_{\bullet}]]$ and torsion over ${A}[[t_{d}]]$ for each $d\in D$. Then
  $\Delta^*\otimes_{{A}[[t_{\bullet}]]}{A}((t_{\bullet}))$ is
  in a natural way an \'{e}tale $(\varphi_{\bullet},\Gamma_{\bullet})$-module over
${A}((t_{\bullet}))$. The functor $\Delta\mapsto
  \Delta^*\otimes_{{A}[[t_{\bullet}]]}{A}((t_{\bullet}))$ is exact.
\end{pro}

{\sc Proof:} For $\Phi(t)=(1+t)^p-1$ see \cite{zabr} Proposition 2.3, which generalizes
a construction of Colmez and Emerton, as recalled in \cite{breuil} Lemma
2.6. We follow the same outline, but we deviate from these texts by
pointing out explicitly the need to invoke a $\psi$-operator on
${A}[[t_{\bullet}]]$. We endow $\Delta^*$ with
an ${A}[[t_{\bullet}]][\Gamma_{\bullet}]$-action by putting $$(a\cdot\ell)(\delta)=\ell(a\delta),$$$$(\gamma\cdot \ell)(\delta)=\ell(\gamma^{-1}\delta)$$for $a\in
{A}[[t_{\bullet}]]$, $\ell\in \Delta^*$, $\delta\in\Delta$ and $\gamma\in\Gamma_{\bullet}$.

{\it Step 1.} Here we assume that the $A$-action on $\Delta$ factors through $k$. Choose a $k$-vector space complement
of ${\Delta}[t_{\bullet}]$ in ${\Delta}$ and let ${\Delta}_0^*$ be the sub
vector space of ${\Delta}^*$ consisting of linear forms vanishing on it. By
the topological Nakayama Lemma, ${\Delta}_0^*$ generates ${\Delta}^*$ as a
$k[[t_{\bullet}]]$-module. Thus, ${\Delta}^*$
is finitely generated as a $k[[t_{\bullet}]]$-module. Similarly,
$(k[[t_{\bullet}]]\otimes_{\varphi_{d},k[[t_{\bullet}]]}\Delta)^*$ is
finitely generated as a $k[[t_{\bullet}]]$-module, and it has the same
generic rank as ${\Delta}^*$. \footnote{For a $k[[t_{\bullet}]]$-module $X$, resp. a $k((t_{\bullet}))$-module $X$, we call ${\rm dim}_{{\rm Frac}(k[[t_{\bullet}]])}X\otimes_{k[[t_{\bullet}]]} {\rm Frac}(k[[t_{\bullet}]])$, resp. ${\rm dim}_{{\rm Frac}(k[[t_{\bullet}]])}X\otimes_{k((t_{\bullet}))} {\rm Frac}(k[[t_{\bullet}]])$, the generic rank of $X$, where ${\rm Frac}(k[[t_{\bullet}]])$ denotes the fraction field of $k[[t_{\bullet}]]$.}

{\it Step 2.} Fix $d\in D$. The map ${A}[[t_{\bullet}]]\otimes_{\varphi_{d},{A}[[t_{\bullet}]]}\Delta\stackrel{{\rm
  id}\otimes \varphi_{d}}{\longrightarrow}\Delta$ gives rise to the ${A}((t_{\bullet}))$-linear
map\begin{gather}\Delta^*\otimes_{{A}[[t_{\bullet}]]}{A}((t_{\bullet}))\stackrel{({\rm
  id}\otimes \varphi_{d})^*\otimes {A}((t_{\bullet}))
}{\longrightarrow}({A}[[t_{\bullet}]]\otimes_{\varphi_{d},{A}[[t_{\bullet}]]}\Delta)^*\otimes_{{A}[[t_{\bullet}]]}{A}((t_{\bullet})).\label{emcobrza1}\end{gather}We claim that it is
bijective. By a standard devissage argument, to do this we may assume that the $A$-action on $\Delta$ factors through $k$, and then the map (\ref{emcobrza1}) reads\begin{gather}\Delta^*\otimes_{k[[t_{\bullet}]]}k((t_{\bullet}))\stackrel{({\rm
  id}\otimes \varphi_{d})^*\otimes k((t_{\bullet}))
}{\longrightarrow}(k[[t_{\bullet}]]\otimes_{\varphi_{d},k[[t_{\bullet}]]}\Delta)^*\otimes_{k[[t_{\bullet}]]}k((t_{\bullet})).\label{emcobrza11}\end{gather}Let $C_{d}$ be the cokernel of
$k[[t_{\bullet}]]\otimes_{\varphi_{d},k[[t_{\bullet}]]}\Delta\stackrel{{\rm
  id}\otimes \varphi_{d}}{\longrightarrow}\Delta$. As $\Delta$ is finitely generated
over $k[[t_{\bullet}]][\varphi_{\bullet}]$ we see that $C_{d}$ is
finitely generated over the subring $k[[t_{\bullet}]][\varphi_{d'}]_{d'\ne
  d}$. With $\Delta$ also $C_{d}$ is a torsion $k[[t_{d}]]\subset
k[[t_{\bullet}]]$-module. As $t_{d}$ belongs to the center of $k[[t_{\bullet}]][\varphi_{d'}]_{d'\ne
  d}$ we obtain that $C_{d}$ is killed by some power of $t_{d}$, hence
$C_{d}^*\otimes_{k[[t_{\bullet}]]}k((t_{\bullet}))=0$. It follows that
the map (\ref{emcobrza11}) is
injective. By what we saw in step 1, the source and the
target of the map (\ref{anneh1}) are finitely generated over
$k((t_{\bullet}))$, and their generic ranks coincide. Thus, the map (\ref{emcobrza11}) has a cokernel which is a
finitely generated torsion module over $k((t_{\bullet}))$. Its annihilator
is a $\Gamma_{\bullet}$-invariant ideal as the map (\ref{emcobrza11}) is
$\Gamma_{\bullet}$-equivariant, hence it must be all of $k((t_{\bullet}))$
by Theorem \ref{gergelytrick}. Therefore the map
(\ref{emcobrza11}) is an isomorphism. 

{\it Step 3.} We use the $\psi$-operator on ${A}[[t_{\bullet}]]$ to define
the ${A}[[t_{\bullet}]]$-linear map\begin{gather}
{A}[[t_{\bullet}]]\otimes_{\varphi_{d},{A}[[t_{\bullet}]]}(\Delta^*)\longrightarrow({A}[[t_{\bullet}]]\otimes_{\varphi_{d},{A}[[t_{\bullet}]]}\Delta)^*,\label{anneh}\\a\otimes\ell\quad\mapsto\quad[b\otimes
  x\mapsto \ell(\psi_{d}(ab)x)].\notag\end{gather}(Note that the action of ${A}[[t_{\bullet}]]$ on the target of the map  (\ref{anneh}) is via the {\it arguments} of the linear forms.) We claim that its base change to ${A}((t_{\bullet}))$ is
bijective. By a standard devissage argument, to do this we may assume that the $A$-action on $\Delta$ factors through $k$, and then the map (\ref{anneh}) reads\begin{gather}k[[t_{\bullet}]]\otimes_{\varphi_{d},k[[t_{\bullet}]]}(\Delta^*)\longrightarrow(k[[t_{\bullet}]]\otimes_{\varphi_{d},k[[t_{\bullet}]]}\Delta)^*.\label{anneh1}\end{gather}

As $k[[t_{\bullet}]]$ is free over
$\varphi_{d}(k[[t_{\bullet}]])$ with basis $1,t_d,\ldots, t_d^{q-1}$, each element in $k[[t_{\bullet}]]\otimes_{\varphi_{d},k[[t_{\bullet}]]}(\Delta^*)$
can be written in the form $\sum_{i=0}^{q-1}t_d^i\otimes\ell_i$ with $\ell_i\in
\Delta^*$. Suppose that $\ell_{i_0}\ne0$ for some $0\le i_0\le q-1$. Choose some
$x\in \Delta$ with $\ell_{i_0}(x)\ne0$. If we are in one of the cases

(i) $F\ne {\mathbb Q}_p$, or

(ii) $F={\mathbb Q}_p$ and $i_0\notin\{0,q-1\}$, or

(iii) $F={\mathbb Q}_p$ and $i_0=q-1$ and $\ell_0=0$,

then the map (\ref{anneh1}) takes $\sum_{i=0}^{q-1}t_d^i\otimes\ell_i$ to a
linear form on
${A}[[t_{\bullet}]]\otimes_{\varphi_{d},{A}[[t_{\bullet}]]}\Delta$ which takes
$t^{q-1-i_0}\otimes x$
to $$\sum_{i=0}^{q-1}\ell_i(\psi_d(t^{q-1-i_0+i})x)=\ell_{i_0}(\psi_d(t^{q-1})x)=\ell_{i_0}(x)\ne0.$$This
follows from formulae (\ref{psineqp}) and (\ref{psieqqp}). If
however we are in the case

(iv) $F={\mathbb Q}_p$ and $i_0=0$ and $\ell_1=0$

then the map (\ref{anneh1}) takes $\sum_{i=0}^{q-1}t_d^i\otimes\ell_i$ to a
linear form on
${A}[[t_{\bullet}]]\otimes_{\varphi_{d},{A}[[t_{\bullet}]]}\Delta$ which takes
$t^{q-1}\otimes x$
to $$\sum_{i=0}^{q-1}\ell_i(\psi_d(t^{q-1+i})x)=\ell_0(\psi_d(t^{q-1})x)=\ell_{0}(x)\ne0,$$again
due to formula (\ref{psieqqp}). Since (choosing $i_0$ suitably if $F={\mathbb Q}_p$) at least one of the cases
(i) --- (iv) must occur, we have shown that the map (\ref{anneh1}) is
injective. By what we saw in step 1, the source and the
target of the map (\ref{anneh1}) are finitely generated over
$k[[t_{\bullet}]]$, and their generic ranks coincide. Thus, the base change
of the map (\ref{anneh1}) to $k((t_{\bullet}))$ has a cokernel which is a
finitely generated torsion module over $k((t_{\bullet}))$. Its annihilator
is a $\Gamma_{\bullet}$-invariant ideal as the map (\ref{anneh1}) is
$\Gamma_{\bullet}$-equivariant, hence it must be all of $k((t_{\bullet}))$
by Theorem \ref{gergelytrick}. Therefore the base change of the map
(\ref{anneh1}) to $k((t_{\bullet}))$ is an isomorphism. 

{\it Step 4.} Composing the base change to ${A}((t_{\bullet}))$ of
the map (\ref{anneh}) with the inverse of the map (\ref{emcobrza11}) we obtain
a bijective ${A}((t_{\bullet}))$-linear
map \begin{gather}{A}((t_{\bullet}))\otimes_{\varphi_{d},{A}((t_{\bullet}))}(\Delta^*\otimes_{{A}[[t_{\bullet}]]}{A}((t_{\bullet})))=({A}[[t_{\bullet}]]\otimes_{\varphi_{d},{A}[[t_{\bullet}]]}(\Delta^*))\otimes_{{A}[[t_{\bullet}]]}{A}((t_{\bullet}))\notag\\\longrightarrow
  \Delta^*\otimes_{{A}[[t_{\bullet}]]}{A}((t_{\bullet})).\notag\end{gather}This
  yields the desired operator
$\varphi_{d}$ on $\Delta^*\otimes_{{A}[[t_{\bullet}]]}{A}((t_{\bullet}))$. It is immediate that $\Delta^*\otimes_{{A}[[t_{\bullet}]]}{A}((t_{\bullet}))$ is
  in a natural way an \'{e}tale $(\varphi_{\bullet},\Gamma_{\bullet})$-module over
${A}((t_{\bullet}))$.

The exactness of $\Delta\mapsto
  \Delta^*\otimes_{{A}[[t_{\bullet}]]}{A}((t_{\bullet}))$ is clear.\hfill$\Box$\\

{\bf Remark:} In practice it may turn out to be helpful if the $t_{d}$ acted
surjectively on $\Delta$, or if at least, for the purpose of constructing the
associated \'{e}tale $(\varphi_{\bullet},\Gamma_{\bullet})$-module in
Proposition \ref{nopsi}, $\Delta$ can be replaced by some $\underline{\Delta}$
which has this property. Let us explain why this is indeed possible at least if $|D|=1$. We omit
indices $(.)_{d}$. In the same way as before, the construction of the $\varphi$-operator on $\Delta^*\otimes_{{A}[[t]]}{A}((t))$ relies on bijectivity claims for certain ${A}((t))$-linear maps, and these claims are reduced to similar claims on
$k((t))$-linear maps. Thus, it is enough to discuss a $k[[t]][\varphi,\Gamma]$-module $\Delta$ which is finitely
  generated over $k[[t]][\varphi]$, admissible and torsion over
  $k[[t]]$. We then let $\underline{\Delta}=\cap_{n\ge0}t^n\Delta$. Notice that
$\underline{\Delta}=\Delta$ if $t$ acts surjectively on $\Delta$.

{\it Claim: $t$ acts surjectively on $\underline{\Delta}$.}

Let $x\in \underline{\Delta}$. For $n\ge0$ choose some $y_n\in {\Delta}$
with $t^ny_n=x$. As ${\Delta}[t]$ is finite dimensional and $k$ is finite,
${\Delta}[t]$ is a finite set, hence also $t^{-1}(x)$ is a finite subset of
$\Delta$. It contains the set of all $t^{n-1}y_n$ with $n\in{\mathbb N}$. Therefore, there
must exist an infinite subset $X$ of ${\mathbb N}$ and some $y\in t^{-1}(x)$
with $y=t^{n'-1}y_{n'}$ for all $n'\in X$. In particular we have $y\in
\underline{\Delta}$. The claim is proven.

The inclusion $\underline{\Delta}\to {\Delta}$ induces a surjection
${\Delta}^*\to \underline{\Delta}^*$ whose cokernel is $t$-torsion, hence
${\Delta}^*\otimes_{k[[t]]}k((t))=\underline{\Delta}^*\otimes_{k[[t]]}k((t))$,
and similarly
$(k[[t]]\otimes_{\varphi,k[[t]]}{\Delta})^*\otimes_{k[[t]]}k((t))\cong
(k[[t]]\otimes_{\varphi,k[[t]]}\underline{\Delta})^*\otimes_{k[[t]]}k((t))$.\\

Taking $D$ with $|D|=1$ and omitting indices $(.)_{d}$ we obtain the $A$-algebra $A((t))[\varphi,\Gamma]$. For general $D$ again consider then the $A$-algebra map ${A}[[t_{\bullet}]]\to {A}[[t]]$ with $t_{d}\mapsto t$ for
each $d\in D$. It induces an ${A}$-algebra map ${A}((t_{\bullet}))\to
{A}((t))$, and given an
${A}((t_{\bullet}))[\varphi_{\bullet},\Gamma_{\bullet}]$-module ${\bf
  D}$ we may regard ${A}((t))\otimes_{{A}((t_{\bullet}))}{\bf
  D}$ as an ${A}((t))[\varphi,\Gamma]$-module by letting $\varphi$
(resp. $\gamma\in\Gamma$) act as $1\otimes\prod_{d\in D}\varphi_{d}$
(resp. as $1\otimes\prod_{d\in D}\gamma_{d}$).

\begin{pro}\label{gergelyred} The functor ${\bf
  D}\mapsto {A}((t))\otimes_{{A}((t_{\bullet}))}{\bf
  D}$ is
  an exact and faithful functor from the category of \'{e}tale $(\varphi_{\bullet},\Gamma_{\bullet})$-modules over ${A}((t_{\bullet}))$ to the category of \'{e}tale $(\varphi,\Gamma)$-modules over ${A}((t))$.
\end{pro}

{\sc Proof:} Again we follow \cite{zabr} Proposition 2.6 (which treats the
case $\Phi(t)=(1+t)^p-1$). For exactness we argue by induction on $|D|$, the case $|D|=1$ being
tautological. If $|D|>1$ pick $d_1\ne d_2$ in $D$. Then $t_{d_1}-t_{d_2}$
lies in the kernel of ${A}((t_{\bullet}))\to
{A}((t))$ and $${A}((\{t_{d}\}_{d\ne d_1}))\longrightarrow
{A}((t_{\bullet}))/(t_{d_1}-t_{d_2})$$is bijective. Let $0\to{\bf
  D}_1\to{\bf
  D}_2\to{\bf
  D}_3\to0$ be an exact sequence of \'{e}tale
$(\varphi_{\bullet},\Gamma_{\bullet})$-modules over
${A}((t_{\bullet}))$. Once we know that $$\frac{{\bf
  D}_1}{(t_{d_1}-t_{d_2}){\bf
  D}_1}\longrightarrow\frac{ {\bf
  D}_2}{(t_{d_1}-t_{d_2}){\bf
  D}_2}$$ is injective we are done, by appealing to the induction
  hypothesis. Let ${\bf D}_3'$ be the ${A}((t_{\bullet}))[\Gamma_{\bullet}]$-submodule of ${\bf D}_3$ generated by the elements $z\in{\bf D}_3$ for which there is some $y\in{\bf D}_2$ with  $(t_{d_1}-t_{d_2})y\in{\bf D}_1$. For all such $z$ we have $(t_{d_1}-t_{d_2})z=0$. Hence, as ${\bf D}_3'$ is a finitely generated ${A}((t_{\bullet}))$-module, there is some $\alpha\in {A}((t_{\bullet}))$ mapping to a non-zero element in ${k}((t_{\bullet}))$ and with $\alpha{\bf D}_3'=0$. Corollary \ref{herzmariens} then says ${\bf D}_3'=0$, and this is what we needed. Similarly, faithfulness follows once more from Theorem \ref{gergelytrick},
just as in \cite{zabr} Proposition 2.8.\hfill$\Box$\\

\begin{flushleft} \textsc{Humboldt-Universit\"at zu Berlin\\Institut f\"ur Mathematik\\Rudower Chaussee 25\\12489 Berlin, Germany}\\ \textit{E-mail address}:
gkloenne@math.hu-berlin.de \end{flushleft} \end{document}